\documentclass[12pt]{amsart}
\usepackage[]{amsmath, amsthm, amsfonts, amssymb, epsfig}
\newcommand\C{{\mathbb C}}

\newcommand\dee{\partial}
\newcommand\Om{\Omega}
\newcommand\Obar{\overline{\Omega}}
\newcommand\Ot{\widetilde\Omega}

\renewcommand\phi{\varphi}

\numberwithin{equation}{section}

\begin{document}

\title[Temperature sensing and Bergman kernel]{
Remote temperature sensing in 2D and the Bergman kernel
}
\author[S.R.~Bell]{Steven R.~Bell}
\author[L.~McNabb]{Leah McNabb}

\address{Steven Bell  \hfill\break\indent
Mathematics Department \hfill\break\indent
Purdue University \hfill\break\indent
West Lafayette, IN 47907}
\email{bell@math.purdue.edu}

\address{Leah McNabb \hfill\break\indent
Mathematics Program\hfill\break\indent
Seton Hill University\hfill\break\indent
Greensburg, PA 15601
}
\email{lmcnabb@setonhill.edu}


\begin{abstract}
We explore the problem of estimating the steady
state temperature in a two-dimensional domain
at a point knowing the temperature to high order
at another point. We find connections to the
Bergman kernel of the domain, Runge's theorem,
and approximate null quadrature identities.
\end{abstract}

\maketitle

\begin{quotation}
\begin{center}
\begin{em}
To celebrate the life and work of Marek Jarnicki
\end{em}
\end{center}
\end{quotation}

\theoremstyle{plain}

\newtheorem {thm}{Theorem}[section]
\newtheorem {lem}[thm]{Lemma}

\hyphenation{bi-hol-o-mor-phic}
\hyphenation{hol-o-mor-phic}

\section{Introduction}
\label{sec1}
Suppose that $\Om$ is a bounded domain in the complex
plane and $u$ is a real valued harmonic function on
$\Om$ (which we can think of as representing the
steady-state temperature in a $\Om$-shaped metal plate).
Suppose further that $u$ is bounded by $M>0$ on $\Om$, i.e.,
$|u(z)|<M$ for all $z\in\Om$, and assume that $a$ and
$b$ are distinct points in $\Om$. We know that if $u$
vanishes to infinite order at $a$, then $u(b)=0$ too.
We will explore the question, given $\epsilon>0$, is
it possible to find a positive integer $N$ and real
coefficients $c_{nm}$ for $0\le n,m\le N$ such that
$$\left|u(b)-\sum_{0\le n,m\le N} c_{nm}
\frac{\dee^{n+m}u}{\dee^n x\,\dee^m y}(a)\right|<\epsilon?$$
We emphasize that $N$ and the coefficients $c_{nm}$
might depend on $a$,$b$, $\epsilon$, and $M$, but not
on $u$.  We like to think of this problem as follows.
Suppose we know that the (steady-state) temperature
in the two-dimensional state of Indiana is bounded by
$130^\circ$. If we can measure the temperature in
Indianapolis to a high enough degree of accuracy to be able to
compute finitely many partial derivatives, can
we estimate the temperature in West Lafayette to
within one-half degree?

We remark that one never needs to use difference quotients to
compute partial derivatives of a harmonic function. One can
always reduce them to an integral of the function times a smooth
function on a circle about the point of interest (or smooth simple
closed curve containing the point) via the Poisson integral
formula.

We will investigate a method to solve the sensing problem
that involves the Bergman kernel and ``approximate null
quadrature identities'' for holomorphic functions that will
lead to some interesting connections to Runge's theorem and
approximation problems for holomorphic functions.

\section{Proof of concept}
\label{sec2}
We will always assume that $\Om$ is a bounded domain in
the complex plane and that $a$ and $b$ are distinct points
in $\Om$. Let $K(z,w)$ denote the Bergman kernel
for $\Om$ and write $K_a(z)=K(z,a)$. Also write $K_a^0(z)=K_a(z)$
and
$$K_a^m(z)=\left.\frac{\dee^m}{\dee\bar w^m}K(z,w)\right|_{w=a}.$$
If $h$ is a function in the Bergman space, then
$$h(z)=\langle h,K_z\rangle,$$
the inner product being the $L^2$ inner product on $\Om$
with respect to area measure. Furthermore,
$$h^{(m)}(z)=\langle h,K_z^m\rangle.$$
Hence, if a function in the Bergman space is orthogonal to
$K_a^m$ for $m=0,1,2,\dots$, then it vanishes to infinite
order at $a$ and is consequently identically zero on $\Om$.
This implies that the complex linear span of the functions
$$\{K_a^m(z)\,:\, m=0,1,2,\dots\}$$
is dense in the Bergman space. Therefore, given $\epsilon>0$,
there is an element in this span that is within $\epsilon$
of $K_b$ in $L^2$ norm. Say $\lambda:=K_b-\sum_{m=0}^N c_m K_a^m$
has $L^2$ norm less than $\epsilon$. Then,
\begin{equation}
\label{eqn1}
\langle h,\lambda\rangle=h(b)-\sum_{m=0}^N {\bar c}_m h^{(m)}(a)
\end{equation}
for functions $h$ in the Bergman space, and we can estimate
\begin{equation}
\label{eqn2}
\left|
h(b)-\sum_{m=0}^N {\bar c}_m h^{(m)}(a)
\right|
\le \epsilon\|h\|
\end{equation}
where the norm is the $L^2$ norm. Consequently, if $|h|<M$ on
$\Om$, then we have bounded the left hand side of this
inequality by $\epsilon M$ times the square root of the area
of $\Om$. We have solved the approximation problem in the
introduction in the realm of holomorphic functions!
This is an abstract and existential proof of concept. We will present a
constructive proof in the next section where the number $N$ and
constants $c_{nm}$ can be determined and related to the parameters.

We like to think of (\ref{eqn1}) and (\ref{eqn2}) as representing an
``approximate null quadrature identity'' on $\Om$.
We will now show how expressions like (\ref{eqn1}) and
(\ref{eqn2}) can be used to
deduce a solution to the harmonic function approximation problem.
The idea is to start with a real valued harmonic function $u$ on
$\Om$, let $v$ be a harmonic conjugate of $u$ on a simply connected
subdomain $\Om_1$ compactly contained in $\Om$ that contains $a$
and $b$ with $v(a)=0$, and then apply an estimate
like (\ref{eqn2}) to the holomorphic function $h:=u+iv$ on $\Om_1$.
The key is an argument motivated by the proof given by the second
author in her Purdue PhD thesis \cite{McN} that one-point quadrature
domains for analytic functions are also one-point quadrature
domains for harmonic functions. Since $u$ is harmonic on
$\Om_1$ and $\Om_1$ is simply connected, the harmonic
conjugate $v$ can be defined on $\Om_1$ as a line integral
$$v(z)=\int_{\gamma_a^z} (-u_y,u_x)\cdot d{\vec s},$$
where $\gamma_a^z$ denotes any smooth curve in $\Om_1$ that starts
at $a$ and ends at $z$. Note that $v(a)=0$, as desired. The proof of
the Cauchy-Riemann equations yields that $h'(z)=u_x+iv_x$
and iterating yields that $h^{(m)}(z)=\frac{\dee^m u}{\dee x^m}
+i\frac{\dee^m v}{\dee x^m}$. Note that the Cauchy-Riemann
equations say that $v_x=-u_y$ and differentiating this equation
$m-1$ times with respect to $x$ yields 
$\frac{\dee^m v}{\dee x^m}=
-\frac{\dee^{m} u}{\dee x^{m-1}\,\dee y}$
for $m>1$. Hence,
$$h' = u_x-iu_y\qquad \text{and} \qquad
h^{(m)} =\frac{\dee^m u}{\dee x^m}
-i\frac{\dee^{m} u}{\dee x^{m-1}\,\dee y}$$if $m>1$.
Assume that an approximate null quadrature identity
of the form (\ref{eqn2}) holds for our function $h$
in the Bergman space associated to $\Om_1$.
Write $c_m=A_m+iB_m$ in terms of real and imaginary
parts.
Now note that the real part of
\begin{align*}
h(b)- & \sum_{m=0}^N {\bar c}_m h^{(m)}(a) = \\
 & h(b)-(A_0-iB_0)u(a)-\sum_{m=1}^N (A_m-iB_m)\left[
\frac{\dee^m u}{\dee x^m}-i
\frac{\dee^{m} u}{\dee x^{m-1}\,\dee y}\right]
\end{align*}
can be expressed as $u(b)$ minus a {\it real\/} linear combination 
of $u$ and some of its partial derivatives of order less than
or equal to $N$ at $a$. To be precise, the real part is
$$u(b)-A_0u(a)-A_1 u_x(a)+B_1u_y(a)
-\sum_{m=2}^N \left[A_m
\frac{\dee^m u}{\dee x^m}
-B_m\frac{\dee^{m} u}{\dee x^{m-1}\,\dee y}\right](a),$$
which is bounded in absolute value by $\epsilon\|h\|$, where the
norm is the $L^2$ norm on $\Om_1$. We now claim that $\|h\|$ is
less than or equal to a positive constant times the $L^2$ norm of $u$
on the larger set $\Om$. Since clearly $|u|\le |h|$ on $\Om_1$, we are
reduced to estimating $|v|$. Estimating the harmonic conjugate of
a given harmonic function on a subset of a domain makes a good
exercise in a graduate course in complex analysis. One approach
could be to use the line integral formula for $v$ to show that we
can bound $v$ on $\Om_1$ by estimating first partial derivatives
of $u$ on a slightly larger relatively compact open set $\Om_2$ of
$\Om$ that contains the closure of $\Om_1$. To do this,
we might cover the closure of $\Om_1$ by finitely many small discs
of radius $d/2$ where $d$ is less than the distance of $\Om_1$ to
the boundary of $\Om_2$. The value of the first partial derivatives of $u$ 
inside the discs of radius $d/2$ can be expressed as integrals
of $u$ against smooth functions on the circles of radius $d$
bounding the discs by differentiating under the integral in the
Poisson formula. This completes the proof that $N$ and appropriate
constants {\it exist\/} in the temperature sensing problem.

We now turn to describing a constructive proof of the holomorphic
sensing problem, which would be the basis of a constructive proof
of the harmonic sensing problem using the techniques of the paragraph
above.

\section{A more constructive approach}
\label{sec3}
The Hilbert space arguments in \S2 are purely existential. In this
section, we describe a method to explicitly determine the order
$N$ and to compute the coefficients in an identity like (\ref{eqn2})
for holomorphic functions. The method is based on the observation
that the basic problem is rather straightforward in the unit disc
when $a=0$ and $b\ne 0$. Indeed, since the Bergman
kernel of the unit disc is $K(z,w)=\pi^{-1} (1-z\,\bar w)^{-2}$,
an easy calculation yields $K_0(z)=\pi^{-1}$ and
$K_0^m(z)=\pi^{-1} (m+1)! z^m$. Since $K_b(z)=
\pi^{-1} (1-z\,\bar b)^{-2}$ is holomorphic on the disc of
radius $1/|b|$, we may approximate $K_b(z)$ uniformly
on any disc $D_r(0)$ about the origin of radius $r$ with
$1<r<|b|^{-1}$ by using the Taylor series approximation
(which is a simple derivative of a geometric series).
Indeed
$$K(z,b)=\pi^{-1}\sum_{n=0}^N (n+1)(z\bar b)^n + E_N(z\bar b)$$
where
$$E_N(w)=\pi^{-1}\frac{(N+2)w^{N+1}-(N+1)w^{N+2}}{(1-w)^2},$$
and we see that the error estimates for $E_N(z\bar b)$ are also
simple and easy to control explicitly. To be precise, we get
$$K(z,b)-\sum_{n=0}^N \frac{\bar b^n}{n!}\, K_0^n(z) = E_N(z\bar b)$$
and the error $E_N(z\bar b)$ goes to zero uniformly in $z$ on a
neighborhood of the closed unit disc because $|b|<1$. This is
where the choice of $N$ will be made. However, there will be some
constants that arise in the next paragraph that we demonstrate
how to compute. Some of the constants depend on potentially high
derivatives of a Riemann map. A numerical analyst bent on
implementing our procedure would have to grapple with finding
efficient methods for computing these derivatives. We admit that
our construction is therefore probably of more theoretical than
practical importance unless a particularly simple and explicit formula
for the Riemann map is known. We discuss one such possibility in
\S\ref{sec5}.

The idea now is to choose a relatively compact
simply connected subdomain $\Ot$ of our given bounded
domain $\Om$ that has smooth real analytic boundary and
contains the given points $a$ and $b$ in $\Om$. The Riemann map
$f:\Ot\to D_1(0)$ such that $f(a)=0$ and $f'(a)>0$ extends
holomorphically past the boundary of $\Ot$. Let $B=f(b)$ and let
$F:=f^{-1}$. Let $K(z,w)$ continue to denote the Bergman kernel
of the unit disc and let $\widetilde K(z,w)$ denote the Bergman
kernel of $\Ot$. We will need the transformation
formula for the Bergman kernels in the form
\begin{equation}
\label{transf}
f'(z)K(f(z),w)=
\widetilde K(z,F(w))\,\overline{F'(w)}.
\end{equation}
Note that setting $w=B$ in this formula yields
$f'(z)K_B(f(z))=
\widetilde K_b(z)\,\overline{F'(B)}$,
and since $F'(B)=1/f'(b)$, we have
\begin{equation}
\label{transf2}
\overline{f'(b)}\,f'(z)K_B(f(z))=
\widetilde K_b(z).
\end{equation}
We can differentiate equation (\ref{transf}) $m$-times with
respect to $\bar w$ and set $w=0$ to also obtain
\begin{equation}
\label{transf3}
f'(z)K_0^m(f(z))=
\frac{\dee^m}{\dee\bar w^m}\left[
\widetilde K(z,F(w))\,\overline{F'(w)}
\right]_{w=0},
\end{equation}
which is an explicit complex linear combination of
elements $\widetilde K_a^n$ with $0\le n\le m$.
Observe that by making
$\Lambda(z):=K_B(z)-\sum_{m=0}^N c_m K_0^m(z)$
sufficiently uniformly small on a disc $D_r(0)$ with $r>1$ where
$r$ is sufficiently close to $1$,
since $f'$ extends holomorphically past the boundary,
$f'(z)\Lambda(f(z))$ can also be made uniformly small on an
open set containing the closure of $\widetilde\Om$. We can
now multiply $f'(z)\Lambda(f(z))$ by $\overline{f'(b)}$ and
use the transformation formulas (\ref{transf2}) and
(\ref{transf3})
to obtain our goal of finding an explicit combination
of the form
$\widetilde K_b(z)-
\sum_{m=0}^N \tilde c_{m} {\widetilde K}_a^m(z)$
that is uniformly small on an open subset of $\Om$
containing the closure of $\Ot$.

\section{Connections to Runge's theorem}
\label{sec4}

We will now show that certain Runge approximation problems involving
rational functions without poles on the exterior of a domain that
tend to zero at infinity are closely connected to
Bergman kernel approximation in the interior. Since inversion by a
mapping of the form $1/(z-a)$ where $a$ is a point in the domain
transforms rational functions without poles on the exterior to
rational functions without poles inside a bounded domain, it will also
be seen that rational approximation problems in a domain are equivalent
to Bergman approximation problems. The connection will be
exhibited via the following theorem from \cite[p.~90]{Be},
which has just been waiting for a proper application.

Suppose $\Om$ is a bounded simply connected domain bounded by
a simple closed $C^\infty$-smooth curve $\gamma$ and $h$ is
a holomorphic function on $\Om$ in $L^2(\Om)$. Define a
holomorphic function $H(w)$ on the outside of $\Om$ via
$$\overline{H(w)}=\iint_{z\in\Om}\frac{h(z)}{\bar z-\bar w}
\ dA_z,$$
where $dA_z$ represents Lebesgue area measure in the $z$
variable. Notice that $H$ and all its derivatives tend to
zero at infinity.

\begin{thm}
\label{thm1}
Given a positive integer $N$,
there is a positive integer $M$ and a constant $C>0$ such that,
if the $C^M$-norm $\|H\|^M$ of $H$ on the outside of $\Om$
is finite, then $h$ is in $C^N(\Om)$ and the $C^N$-norm
$\|h\|_N$ satisfies
$$\|h\|_N \le C \|H\|^M.$$
Consequently,
$h(z)$ extends $C^\infty$-smoothly to the boundary of $\Om$
from the inside if and only if
$H(w)$ extends $C^\infty$-smoothly to the boundary of $\Om$
from the outside.
\end{thm}

This theorem gives us a direct connection between the problem
of finding an approximate quadrature identity of the form
(\ref{eqn2}) and the problem of approximating $(z-b)^{-1}$
by polynomials in $(z-a)^{-1}$ on the exterior of $\Om$.
Indeed, if $\Om$ is a simply connected domain with smooth boundary
that contains $a$, then $h^{(m)}(a)=\left< h, K_a^m\right>$
when $h$ is in the Bergman space associated with $\Om$.
In particular, if we let $h(z)=(w-z)^{-1}$ where $w$ is a
point outside of the closure of $\Om$, then
$$\frac{m!}{(w-a)^{m+1}}=\left< h, K_a^m\right>$$
and the dual nature of our two problems via Theorem~\ref{thm1}
is clear.

A fascinating possibility is that the way we have suggested
obtaining an approximate quadrature identity using the
Bergman kernel might correspond to a superior way of
doing the classical Runge approximation of $(z-b)^{-1}$ by
polynomials in $(z-a)^{-1}$ on the exterior of $\Om$.

Runge's original proof of one of his famous theorems used an
explicit technique to approximate $1/(z-b)$ by a polynomial
in $1/(z-a)$ when $b\ne a$ on a large open set away from a curve
$\gamma$ connecting $b$ to $a$. The argument involved gradually
moving a point along $\gamma$ to get uniform approximations to
$1/(z-b)$ on the open set given by the set of points that are
more than a distance of $\delta>0$ from the curve. The classic
fact is explained nicely in the textbook \cite[pp.~63-64]{SS}.
We summarize it here. Let $\gamma$ be a curve connecting $b$ to $a$
parametrized by $z(t)$, $\beta\le t\le\alpha$ and let $\delta>0$
be less than $|a-b|$. Let
$U_\delta$ denote the (open) set of points that are a distance
greater than $\delta$ away from the curve $\gamma$. We now
describe step one of the construction. Since $\delta$ is less
than the distance from $b$ to $a$, there is a point $a_1=z(t_1)$ on
$\gamma$ that is a distance $\delta$ from $b$ with $t_1$ as
small as possible. We may approximate the function $(z-b)^{-1}$
by finitely many terms of a Laurent expansion centered at $a_1$
uniformly, say within $\epsilon>0$, outside a disc of radius
$2\delta$ about $a_1$. Let $P_1(z)$ denote the polynomial in
$(z-a_1)^{-1}$ that represents the finitely many terms in the
Laurent expansion. The next step is to replace $b$ by $a_1$,
$\beta$ by $t_1$, and the function $(z-b)^{-1}$ by $P_1(z)$
in step one to obtain a point $a_2$ along the curve and a
polynomial $P_2(z)$ in $(z-a_2)^{-1}$ that is uniformly close,
within $\epsilon$, to $P_1(z)$ outside a disc of radius $2\delta$
about $a_2$. Note that $P_2(z)$ is within $2\epsilon$ of
$(z-b)^{-1}$ outside the intersection of the previous sets
on which the approximations were made. This process can be
repeated until the distance from the $N$-th point $a_N$ to $a$
is less than or equal to $\delta$. In the last step, we let
$a=a_{N+1}$ and approximate $P_N(z)$ outside the disc of radius
$2\delta$ about $a$ by a polynomial in $(z-a)^{-1}$. We can now
replace $\epsilon$ by $\epsilon/(N+1)$ to obtain a uniform approximation
within $\epsilon$ of $(z-b)^{-1}$ by a polynomial in $(z-a)^{-1}$
in the open set $U_{2\delta}$.

We remark that this idea of Runge could have been the basis for a
solution to the remote sensing problem. Indeed if $\Om$ is a simply connected
domain bounded by a smooth curve which compactly contains the
set of points that are within a distance $\delta$ of the curve
$\gamma$, and if the modulus of
$$\frac{1}{z-b}-\sum_{j=1}^N \frac{A_j}{(z-a)^j}$$
is bounded by $\epsilon$ on $U_\delta$, then integrating this
function against a holomorphic function $h\,dz$ around the
boundary of $\Om$ and using the higher order Cauchy integral formulas
yields an approximate one-point quadrature
identity for $h$ that can be bounded by
a constant times $\epsilon$ times the length of the boundary
curve times the supremum of $|h|$ on the boundary. Then the
same ideas used in \S\ref{sec2} to get a solution to the
sensing problem for harmonic functions can be applied to
the holomorphic identity to solve the harmonic problem.

Runge's approximation method is rather crude and likely not
an efficient way to solve the problem, whereas the approximation
scheme on the unit disc via Taylor approximations to a rational
function with poles outside the unit circle is more efficient.
The method on the unit disc pulls back to a method on a simply
connected domain via the Riemann map. It seems possible that
these ideas might find connections to the work of L.~N.~Trefethen
(\cite{Tr,Tr2}) on rational approximation and numerical
analytic continuation.

\section{Further musings on the subject}
\label{sec5}

The techniques of \S\ref{sec2}-\ref{sec3} are not sensitive to
the original region being simply connected. Indeed, the first step was
to shrink the domain down down to a simply connected domain
that contains both points $a$ and $b$. An intriguing idea is
to design ``probe domains'' tailored to solve the temperature
sensing problem where the Riemann maps involved are explicit
and perhaps easier to compute with. Given a possibly unbounded and
non-simply connected domain $\Om$ and points $a$ and $b$ in $\Om$,
let $\gamma$ be a smooth curve in $\Om$ such that $a$
and $b$ are both points interior to the curve, i.e., in the curve,
but not at either endpoint. Next, if $\gamma$ is parametrized by
a smooth complex function $z(t)$, $0\le t\le 1$, approximate the
real and imaginary parts of $z(t)$ by polynomials and make minor
adjustments to get a complex valued polynomial $P(t)$ that agrees
with $z(t)$ at the two image points $a$ and $b$ that stays close to
the original curve and stays in $\Om$.
We may also assume that $P'(t)$ is nonvanishing on $[0,1]$. Let
$\Gamma$ denote the curve in $\Om$ parametrized by $P(t)$. We next
let $t$ go complex so that $P(\tau)$ maps a (perhaps thin) open
simply connected set $U$ about $[0,1]$ in the complex plane
one-to-one onto a snake-like domain about the curve $\Gamma$ that
is compactly contained in $\Om$. Let $E$ be an ellipse (or rectangle,
if preferred) that contains $[0,1]$, is symmetric about the $x$-axis,
and is a subset of $U$. The conformal map from the unit disc onto $E$ is
well known. By composing next with $P(z)$  we get a conformal map
from the unit disc to a ``probe domain'' in $\Om$ that contains
$a$ and $b$. It is easy to arrange that this map sends the origin
to $a$. We may now use the techniques of \S\ref{sec3} to solve
the holomorphic and the heat sensing problems.

It is interesting to note that, in the realm of bounded domains
$\Om$ with $C^\infty$-smooth boundary, it is possible to approximate
$K_b$ by a finite complex linear combination of the functions
$K_a^m$ {\it in the topology of\/} $C^\infty(\Obar)$ (see \cite{Be1}).
Coming up with better constructive ways to do this might lead to
better constructive ways to solve the sensing problems.

It is reasonable to consider the problems of this paper in several
complex variables. Indeed, the topic is closely connected to the
subject of ``one-point quadrature domains,'' and it has been
shown that many of the results in one variable on these domains
generalize to several complex variables (see \cite{Be2} and
\cite{Le}). The soft analysis in the ``proof of concept''
section \S\ref{sec2} yields that there {\it exists\/} a solution
to the holomorphic sensing problem in several variables. The second
author's proof that one-point quadrature domains in $\C$ are
also one point quadrature domains for harmonic functions
carries over to several complex variables to show that
one-point quadrature domains for holomorphic functions in
several variables are also one point quadrature domains for
{\it pluriharmonic\/} functions.  Following the same reasoning
we did in \S\ref{sec2} and restricting to a simply connected
subregion containing the two points yields that we can solve
the sensing problem for {\it pluriharmonic\/}
functions (which have nothing to do with temperature problems).
The fact that there is no Riemann mapping theorem in several
complex variables means alternative means to a more constructive
approximation problem would have to be found.

\end{document}